\input amstex
\documentstyle {amsppt}
\magnification=1200
\vsize=9.5truein
\hsize=6.5truein
\nopagenumbers 
\nologo

\def\cc{\Cal C_{\ge Z}}

\topmatter

\title
Spaces~of~closed~subgroups
of~locally~compact~groups
\endtitle

\author
Pierre de la Harpe
\endauthor

\address
Pierre de la Harpe, Section de Math\'ematiques,
Universit\'e de Gen\`eve, C.P. 64,
\newline
CH--1211 Gen\`eve 4.
%  \newline
E--mail: Pierre.delaHarpe\@math.unige.ch
\endaddress

\date
16th of June, 2008 
% !!! corrections de Courtray faites le 7 juin !!! et 15 juin !!!
\enddate

\leftheadtext
{Spaces of closed subgroups of locally compact groups}
\rightheadtext
{Spaces of closed subgroups of locally compact groups}

\abstract
The set $\Cal C(G)$ of closed subgroups of a locally compact group $G$
has a natural topology which makes it a compact space.
This topology has been defined in various contexts
by Vietoris, Chabauty, Fell, Thurston, Gromov, Grigorchuk, 
and many others.

The purpose of the talk was to describe the space 
$\Cal C(G)$ first for a few elementary examples,
then for $G$ the complex plane, in which case $\Cal C(G)$ is a $4$--sphere
(a result of Hubbard and Pourezza),
and finally for the $3$--dimensional Heisenberg group $H$,
in which case $\Cal C(H)$ is a $6$--dimensional singular space
recently investigated by Martin~Bridson, Victor~Kleptsyn and the author
\cite{BrHK}. 

These are slightly expanded notes prepared
for a talk given at several places: the Kortrijk workshop on
{\it Discrete Groups and Geometric Structures, with Applications III,}
May 26--30, 2008; 
the {\it Tripode 14,} 
\'Ecole Normale Sup\'erieure de Lyon, June 13, 2008;
and seminars at the EPFL, Lausanne, and in the 
Universit\'e de Rennes 1.
The notes do not contain any other result than those in \cite{BrHK},
and are not intended for publication.
\endabstract



\endtopmatter

\document

\head{\bf
I.~Mahler (1946)
}\endhead

Let $n$ be a positive integer.
A {\it lattice} in $\bold R^n$ is a subgroup of $\bold R^n$ generated by a basis.
Two lattices $L,L' \subset \bold R^n$ ar {\it close to each other} if there exist
basis $\{e_1,\hdots,e_n\} \subset L$, $\{e'_1,\hdots,e'_n\} \subset L'$
with $e_j$ and $e'_j$ close to each other for each $j$; 
this defines a topology on the space $\Cal L (\bold R^n)$ of lattices in $\bold R^n$.
It coincides with the natural topology on $\Cal L (\bold R^n)$ viewed as
the homogeneous space $GL_n(\bold R) / GL_n(\bold Z)$.
It is easy to check that the covolume function
$L \longmapsto \operatorname{vol}(\bold R^n / L)$ 
and the minimal norm
$L \longmapsto \min L \Doteq \min_{x \in L, x \ne 0} \Vert x \Vert ^2$
are continuous functions $\Cal L (\bold R^n) \longrightarrow \bold R_+^*$.

\proclaim{1.~Mahler's Criterion \cite{Mahl--46}} 
For a subset $\Cal M$ of $\Cal L (\bold R^n)$,
the two following properties are equivalent:
\par
(i) $\Cal M$ is relatively compact;
\par
(ii) there exist 
%two constants 
$C,c > 0$ such that
$\operatorname{vol}(\bold R^n / L) \le C$ and $\min L \ge c$
for all $L \in \Cal M$.
\endproclaim

For proofs of this criterion, see for example \cite{Bore--69, corollaire 1.9}
or \cite{Ragh--72, Corollary 10.9}.

An immediate use of the criterion is the proof that,
in any dimension $n \ge 1$, 
the space $\Cal L^{\text{unimod}}(\bold R^n)$
of lattices of covolume $1$ contains a lattice $L_{\max}$ such that 
$\min L_{\max} = \sup \{ \min L \hskip.1cm \vert \hskip.1cm
L \in \Cal L^{\text{unimod}} (\bold R^n) \}$;
equivalently, there exists a lattice $L_{\max}$
with a maximal density for the associated ball packing.
(It is a much more difficult problem to identify such a lattice $L_{\max}$,
and this problem is open unless $n\le8$ or $n=24$.)

Let us quote one other application of Mahler's Criterion.
It is one of the ingredients of the proof of the following fact, 
special case for $SL_n(\bold R)$ of a more general result
conjectured by Godement and proved independently
by Borel--Harish Chandra and by Mostow--Tamagawa
\cite{Ragh--72, Theorem 10.18}:
{\it an arithmetic lattice in
$SL_n(\bold R)$  is cocompact
if and only if it does not contain any unipotent matrix.}

\bigskip
\head{\bf
II.~Chabauty (1950) and Fell (1962)
}\endhead

For a topological space $X$, 
let $2^X$ denote the set of closed subsets of $X$.
For a compact subset $K$ and a nonempty open subset $U$ of $X$, set
$$
\Cal O_K \, = \, \{ F \in 2^X \hskip.1cm \vert \hskip.1cm F \cap K = \emptyset \}
\hskip.5cm \text{and} \hskip.5cm
\Cal O'_U \, = \, \{ F \in 2^X \hskip.1cm \vert \hskip.1cm F \cap U \ne \emptyset \} .
$$
The finite intersections 
$\Cal O_{K_1} \cap \cdots \cap \Cal O_{K_m} \cap 
\Cal O'_{U_1} \cap \cdots \cap \Cal O'_{U_n}$, $m,n \ge 0$, 
constitute a basis of the {\it Chabauty topology} on $2^X$;
observe that 
$\Cal O_{K_1} \cap \cdots \cap \Cal O_{K_m} = \Cal O_{K_1 \cup \cdots \cup K_m}$.
\par
For example, consider the case $X=\bold R^n$,
and two lattices $L_0,L$ in this space.
If $K$ is a compact subset of $\bold R^n$ disjoint from $L_0$,
then $L \in \Cal O_K$  
if and only if $L$ is also disjoint from $K$;
if $U$ is an open neighbourhood of some point of $L_0$,
then $L \in \Cal O'_U$  
if and only if $L$ has also one point (at least) in $U$.
It follows that the Chabauty topology on $2^{\bold R^n}$
induces on $\Cal L (\bold R^n)$ the same topology
as that considered in Section~I.

\medskip

\Refs\nofrills{}
This topology was defined in \cite{Chab--50} (written in French), 
as a tool to show a version of Mahler's Criterion
which is valid in a large class of topological groups
(Proposition~3 below).
Later, this topology was studied in greater detail by Fell for unrelated purposes, 
precisely for the study of the appropriate topology
on the space of irreducible unitary representations
of a locally compact group; 
Fell (who writes in English) does not quote Chabauty 
(see \cite{Fell--62} and \cite{Fell--64}).
This {\it Chabauty topology,} or {\it Chabauty--Fell topology,} 
or {\it $H$--topology} (terminology of Fell), or 
{\it geometric topology} (terminology of Thurston), on $2^X$
should not be confused 
with other standard topologies on the same space,
of which the study goes back to Hausdorff and Vietoris,
and for which a canonical reference is \cite{Mich--51}
(however, for $X$ compact, the Chabauty topology coincides
with the F--topology of Michael).
The Chabauty topology is also useful 
in the study of low--dimensional manifolds;
see \cite{Thurs, Definition 9.1.1}, 
as well as \cite{CaEG--87, Chapter I.3},
where Chabauty is quoted (but Fell is not).
A pleasant account of the most basic properties of this topology
can be found in \cite{Paul--07}.
\endRefs

\medskip

\proclaim{2.~Proposition} Let $X$ be a topological space,
and let $2^X$ be endowed with the Chabauty topology.
\roster
\item"(i)" 
The space $2^X$ is compact.
\item"(ii)"
In case $X$ is discrete, $2^X$ is  homeomorphic
to the product space $\{0,1\}^X$ with the compact Tychonoff product topology.
\endroster
Suppose moreover that $X$ is a locally compact metric space.
\roster
\item"(iii)"
The Chabauty topology on $2^X$ is induced by the metric defined by
$$
d(F_1,F_2) \, = \, \inf \left\{ \epsilon > 0 \hskip.2cm \Bigg\vert \hskip.2cm
\aligned
&F_1 \cup \Big( X \smallsetminus B(\ast,1/\epsilon) \Big)
\hskip.2cm \subset \hskip.2cm
\Cal V_{\epsilon} \left(  
F_2 \cup \Big( X \smallsetminus B(\ast,1/\epsilon) \Big)
\right)
\\
&F_2 \cup \Big( X \smallsetminus B(\ast,1/\epsilon) \Big)
\hskip.2cm \subset \hskip.2cm
\Cal V_{\epsilon} \left(  
F_1 \cup \Big( X \smallsetminus B(\ast,1/\epsilon) \Big)
\right)
\endaligned
\right\} .
$$
Here, for a subset $S$ of $X$, we write $\Cal V_{\epsilon}(S)$
for $\{x \in X \hskip.1cm \vert \hskip.1cm d(x,S) < \epsilon \}$;
and $B(\ast,1/\epsilon)$ is the open ball of radius $1/\epsilon$
around an arbitrarily chosen base point $\ast \in X$.
\item"(iv)"
A sequence $(F_j)_{j \ge 1}$ of closed subsets of $X$ 
converges in $2^X$ to a closed subset $F$ if and only if
the two following conditions hold:
\item"" \hskip.2cm --- \hskip.2cm
for all $x \in F$, there exists for all $i \ge 1$ a point $x_i \in F_i$
such that $x_i \to x$,
\item"" \hskip.2cm --- \hskip.2cm
for all strictly increasing sequence $(i_j)_{j \ge 1}$
and for all sequences 
\item"" \hskip.8cm
$(x_{i_j})_{j \ge 1}$ such that $x_{i_j} \in F_{i_j}$
and $x_{i_j} \to x \in X$, we have $x \in F$.
\item"(v)"
If $X' = X \sqcup \{\omega\}$ is the one--point compactification of $X$,
then $F \longmapsto F \sqcup \{\omega\}$
is a homeomorphism from $2^X$ to
the subspace of $2^{X'}$ of closed sets containing~$\{\omega\}$.
Moreover, the subset
$\left\{ \{x\} \in 2^X \hskip.1cm \vert \hskip.1cm x \in X \right\} \cup \emptyset$
of $2^X$ and its image in $2^{X'}$ are homeomorphic to $X'$.
\endroster
Suppose moreover that $X=G$ is a locally compact group,
and let $\Cal C (G)$ be the subset of $2^G$ of closed subgroups.
\roster
\item"(vi)"
The subspace $\Cal C (G)$ of $2^G$ is closed,
and therefore compact.
\item"(vii)"
In $\Cal C(G)$, a basis of neighbourhoods of a closed subgroup $C$ is given by
$$
\Cal N_{K,U}(C) \, = \,
\{ D \in \Cal C (G) \mid
D \cap K \subset CU
\hskip.3cm \text{and} \hskip.3cm
C \cap K \subset DU \} .
$$
\endroster
\endproclaim

% \demo{Needs some more thought}
% Space of nonempty closed subsets versus space of closed subsets,
% and so on.
% \enddemo

\Refs\nofrills{}
\demo{Comments}
(i) The proof appears in \cite{Fell--62, Theorem 1}.
In general, $2^X$ need not be Hausdorff,
even if $X$ is metrisable; thus in a french--like terminology,
$2^X$ is quasi--compact.
However, if $X$ is locally compact and possibly non--Hausdorff,
then $2^X$ is a compact Hausdorff space.
In the relevant context and in terms of the geometric topology
(defined by the conditions of (iv)), 
the proof of (i) appears
also in \cite{Thurs, Proposition 9.1.6}, where it is established
that $\Cal C (G)$ is compact for a Lie group $G$.

(ii) For a positive integer $k$, let $F_k$ denote the
free group on $k$ generators. 
The space $\Cal N (F_k)$ of normal subgroups of $F_k$
is closed in $2^{F_k} = \{0,1\}^{F_k}$, and therefore compact.
This space can be naturally identifed with the
{\it space of marked groups on $k$ generators},
namely of groups $\Gamma$ given together with 
a generating set $\{s_1,\hdots,s_k\}$, 
or equivalently together with 
a quotient homomorphism $F_k \longrightarrow \Gamma$.
This space has been intensively studied in recent years:
see among others
\cite{Grom--81, final remarks},
\cite{Grig--84},
\cite{Cham--00},
\cite{ChGu--05},
and \cite{CoGP--07}.

(iii) and (iv) See for example \cite{Paul--07}, Proposition 1.8, Page 60.

(v) See \cite{Fell--62, Page 475},
or \cite{Bour--63}, chapitre VIII, \S~5, exercice~1.

(vi) See \cite{Fell--62, Page 474}.
\enddemo
\endRefs

From now on, $\Cal C (G)$ will denote the compact space of closed subgroups
of a locally compact group $G$, furnished with the Chabauty topology. 
It has several  subspaces of interest, including:
\roster
\item""
the space $\Cal D (G)$ of discrete subgroups of $G$,
\item""
the space $\Cal L (G)$ of lattices of $G$,
\item""
the space $\Cal A (G)$ of closed abelian subgroups of $G$,
\item""
the space $\Cal N (G)$ of closed normal subgroups of $G$.
\endroster
(Recall that a {\it lattice} in $G$
is a discrete subgroup $\Lambda$ such that $G/\Lambda$ has 
a $G$-invariant probability measure.)

\proclaim{3.~Proposition (Chabauty's Mahler's Criterion)}
Let $G$ be a unimodular 
\footnote{
Recall that, if $G$ was not unimodular,
it would not contain any lattice at all;
in other terms, one would have $\Cal L (G) = \emptyset$.
Even if $G$ is unimodular, it may happen that
$\Cal L (G) = \emptyset$; this happens for example
with nilpotent Lie groups of which the Lie algebra $\eufm g$
has no rational form; a {\it rational form} of a Lie algebra $\eufm g$
is a Lie algebra $\eufm g_0$ over $\bold Q$
such that $\eufm g_0 \otimes_{\bold Q} \bold R$
and $\eufm g$ are isomorphic real Lie algebras.
}
locally compact group 
satisfying some extra technical conditions,
for example let $G$ be a connected unimodular Lie group,
and let $M$ be a subset of  $\Cal L (G)$.
Then $\Cal M$ is relatively compact if and only if
\roster
\item"(i)" 
there exists a constant $C > 0$ such that 
$\operatorname{vol}(G/\Lambda) \le C$ for all $\Lambda \in \Cal M$,
\item"(ii)"
there exists a neighbourhood $U$ of $e$ in $G$ such that
$\Lambda \cap U = \{e\}$ for all $\Lambda \in \Cal M$.
\endroster
\endproclaim

   The Chabauty topology provides some natural compactifications.
More precisely, let $S$ be a space of which the points are in natural correspondance
with closed subgroups of a given locally compact group $G$,
in such a way that the corresponding injection 
$\varphi : S \longrightarrow \Cal C (G)$ is continuous.
Then $\overline{\varphi(S)}$ is a  compactification of $S$.
Examples to which this applies are:
\roster
\item"$\to$"
Riemannian symmetric spaces of the non--compact type $S$,
for which each point of $S$ corresponds to a maximal compact subgroup 
of the isometry group  $S$;
\item"$\to$"
Bruhat--Tits buildings;
\item"$\to$"
the space of complete Riemannian manifolds of dimension $n$ 
and of constant sectional curvature $-1$, 
given together with a base point 
and an orthonormal basis of the tangent space at this base point; 
this is a space which can be identified with a space of discrete subgroups
of the isometry group of the hyperbolic space $H^n$.
% (see \cite{Thurs, Chapter~9}).
\endroster

\bigskip
\head{\bf
III.~First examples
}\endhead

If $G = \bold R$, the space $\Cal C (\bold R)$ is homeomorphic 
to a compact interval $[0,\infty]$. 
The points $0$,~$\lambda$ (with $0 < \lambda < \infty$), and $\infty$
correspond respectively to
 the subgroups $\{0\}$, $\frac{1}{\lambda}\bold Z$, and $\bold R$.
 
 \medskip
 
The space $\Cal C (\bold Z)$ is homeomorphic to the subspace
$\left\{\frac{1}{n}\right\}_{n \ge 1} \cup \{0\}$ of $[0,1]$,
with $\frac{1}{n}$ corresponding to $n\bold Z$ and $0$ to $\{0\}$.
[Exercise: the spaces $\Cal C (\bold R / \bold Z)$ and $\Cal C (\bold Z)$
are homeomorphic.]

\medskip

Even if the list of easily understandable spaces $\Cal C (G)$
could be slightly extended 
(exercise: look at $SO(3)$ and at the affine group of the real line), 
it is essentially a very short list.
This was our main motivation 
to understand two more cases~:
the additive group of $\bold C$ ($= \bold R^2$), 
see  \cite{HuPo--79} and Section IV below,
and the Heisenberg group $H$, see \cite{BrHK} and Sections V to VII.

\medskip

In $\bold R^n$, $n \ge 1$, any closed subgroup is isomorphic to one of
$\bold R^a \oplus \bold Z^b$, where $a,b$ are non--negative integers
such that $0 \le a+b \le n$. For a given pair $(a,b)$,
the subspace 
\footnote{
Let $G,A,\hdots,B$ be topological groups.
We denote by $\Cal C_{A,\hdots,B}(G)$ the subspace of $\Cal C (G)$
of closed subgroups of $G$ topologically isomorphic to one of $A,\hdots,B$.
}
$\Cal C_{\bold R^a \oplus \bold Z^b}(\bold R^n)$
of $\Cal C (\bold R^n)$ 
% of groups isomorphic to $\bold R^a \oplus \bold Z^b$
is a homogeneous space of $GL_n(\bold R)$,  for example 
$\Cal C_{\bold Z^n}(\bold R^n) = \Cal L (\bold R^n) = GL_n(\bold R) / GL_n(\bold Z)$,
as already observed in Section~I.
But the way these different \lq\lq strata\rq\rq \ are glued to each other
to form $\Cal C (\bold R^n)$ is complicated, and we do not know of any
helpful description if $n \ge 3$.
\par
To  describe  the case $n=2$, it will be convenient
to think of $\bold C$ rather than of $\bold R^2$.

\bigskip
\head{\bf
IV.~Hubbard and Pourezza (1979)
}\endhead

We will describe $\Cal C (\bold C)$ in three steps:
first the space 
$\Cal C_{\text{nl}} (\bold C) = 
\Cal C_{\{0\},\bold Z,\bold R,\bold R \oplus \bold Z,\bold C}(\bold C)$
of closed subgroups of $\bold C$ which are {\it not} lattices
(easy step), 
then the space $\Cal L (\bold C)$ of lattices
(classical step in complex analysis),
and finally the way these two parts are glued to each other
to form a $4$--sphere (the contribution of Hubbard and Pourezza).

\medskip
\subhead
IV.1.~The space $\Cal C_{\text{nl}} (\bold C)$ is a $2$--sphere
\endsubhead

In $\Cal C (\bold C)$, the subset of closed subgroups isomorphic to $\bold R$
is the real projective line $\bold P^1$, namely a circle.
Each infinite cyclic 
% closed 
subgroup of $\bold C$
is contained in a unique real line, 
and the subgroup $\{0\}$ is contained in all of them.
It follows that the space $\Cal C_{ \{0\}, \bold Z, \bold R }(\bold C)$
is a cone over $\bold P^1$, namely a closed disc.
We rather think of it as the closed lower hemisphere 
of a $2$--sphere $\bold S^2$.

A closed subgroup $C$ of $\bold C$ isomorphic to $\bold R \oplus \bold Z$
is uniquely determined by
its connected component $C^o$, isomorphic to $\bold R$,
and by the minimal norm $\min_{z \in C, z \notin C^o} \vert z \vert$.
When this minimal norm is very large [respectively very small],
$C$ is \lq\lq near\rq\rq \ a closed subgroup of $\bold C$
isomorphic to $\bold R$ [respectively is \lq\lq near\rq\rq \ $\bold C$].
It follows that the space $\Cal C_{\bold R, \bold R \oplus \bold Z, \bold C}(\bold C)$
is also a cone over $\bold P^1$, which can be identified
to the closed upper  hemisphere of $\bold S^2$.

Consequently, $\Cal C_{\text{nl}}(\bold C)$ is homeomorphic to $\bold S^2$, with
\roster
\item""
$\{0\}$ corresponding to the South Pole,
\item""
$\Cal C_{\bold Z}(\bold C)$ to the complement of the South Pole
in the open lower hemisphere,
\item""
$\Cal C_{\bold R}(\bold C)$ to the equator,
\item""
$\Cal C_{\bold R \oplus \bold Z}(\bold C)$ to the complement of the North Pole
in the open upper hemisphere,
\item""
$\bold C$ to the North Pole.
\endroster

\medskip
\subhead
IV.2.~The space $\Cal L (\bold C)$ is the product of an open interval
with the complement of a trefoil knot in $\bold S^3$
\endsubhead

For a lattice $L$ in $\bold C$, set as usual
$$
g_2(L) \, = \, 60 \sum_{z \in L, z \ne 0} z^{-4} ,
\hskip.5cm 
% \text{and} \hskip.5cm
g_3(L) \, = \, 140 \sum_{z \in L, z \ne 0} z^{-6} ,
\hskip.5cm
\Delta(L) \, = \, g_2(L)^3 - 27 g_3(L)^3.
$$
% and $\Delta(L) \, = \, (g_2(L))^3 - 27 (g_3(L))^2$.
Let $\Sigma$ be the complex curve in $\bold C^2$ of equation $a^3 = 27 b^2$.
It is a classical fact that we have a homeomorphism
$$
g \, : \,
\left\{ \aligned
\Cal L (\bold C) \hskip.2cm &\longrightarrow \hskip.8cm \bold C^2 \smallsetminus \Sigma
\\
L \hskip.5cm &\longmapsto \hskip.2cm (g_2(L),g_3(L)) 
\endaligned \right.
$$
(see for example \cite{SaZy--65, \S~VIII.13}).
By the same formulas, $g$ can be extended to a homeomorphism
$$
g' \, : \, \Cal L_{\{0\},\bold Z,\bold Z^2}(\bold C) 
\hskip.2cm \longrightarrow \hskip.2cm
\bold C^2 .
$$
Observe that $\Sigma$, viewed as a real surface in $\bold C^2 \approx \bold R^4$,
is smooth outside the origin, and that its intersection with the unit sphere $\bold S^3$
of equation $\vert a \vert^2 + \vert b \vert^2 = 1$ is a {\it trefoil knot}
$$
T \, = \, \left\{ (a,b) \in \bold S^3 \hskip.1cm \vert \hskip.1cm a^3 = b^2 \right\} .
$$
(It follows that the origin is indeed a singular point of $\Sigma$.)
\par

The multiplicative group $\bold C^*$ acts on $\Cal C (\bold C)$ and on $\bold C^2$
by
$$
(s,C) \hskip.2cm \longmapsto \hskip.2cm \sqrt{s}C
\hskip.5cm \text{and} \hskip.5cm
\big( s,(a,b)\big) \hskip.2cm \longmapsto \hskip.2cm (s^{-2}a,s^{-3}b) 
$$
(observe that $\sqrt{s}C$ is well defined, because $-C=C$).
Moreover, the homeomorphism $g'$ is $\bold C^*$--equivariant.
\par

Let us restrict these actions  to the subgroup $\bold R^*_+$ of $\bold C^*$.
The resulting action of $\bold R^*_+$ on $\Cal L (\bold C)$
is free, and its  orbits are transverse to the subset
$\Cal L^{\text{unimod}}(\bold C)$ of unimodular lattices in $\bold C$.
The resulting action of $\bold R^*_+$ on $\bold C^2 \smallsetminus \Sigma$
is also free, and its orbits are transverse to the sphere~$\bold S^3$.
And the homeomorphism $g$ is $\bold R^*_+$--equivariant.
[The restriction of the action to the subgroup of complex numbers
of modulus one is also interesting, producing on $\Cal L^{\text{unimod}}(\bold C)$
the structure of a Seifert manifold, but we will not expand this here.]
\par
However, note that the transversals $\Cal L^{\text{unimod}}(\bold C)$
and $\bold S^3 \smallsetminus T$ {\it do not} 
correspond to each other by the homeomorphism $g$.
Compare with the first comment following Theorem~4 below.

% SUITE DEJA DITE PLUS BAS
% \medskip
% 
% \Refs\nofrills{}
% It is true, even if not used below, 
% that $\Cal L^{\text{unimod}}(\bold C)$ is a homogeneous space
% of the universal covering group of $SL_2(\bold R)$ by a discrete subgroup,
% and is consequently an aspherical space.
% The space $\Cal L (\bold C)$, which is homeomorphic to the direct product
% of the transversal $\Cal L^{\text{unimod}}(\bold C)$ and a fiber $\bold R^*_+$,
% is also aspherical.
% \endRefs

\medskip
\subhead
IV.3.~The space $\Cal C (\bold C)$ is a $4$--sphere
\endsubhead

The covolume is usually defined on the set of lattices, 
but there is no difficulty to see it as a continuous mapping
$\Cal C_{\{0\},\bold Z,\bold Z^2}(\bold C) \longrightarrow [0,\infty]$,
cyclic subgroups being of infinite covolume.
We can therefore consider the set
$\Cal C^{\text{covol} \ge 1}(\bold C)$
of closed subgroups of $\bold C$ of covolume at least $1$;
it is a subspace of $\Cal C_{\{0\},\bold Z,\bold Z^2}(\bold C)$ which contains
$\Cal C_{\{0\},\bold Z}(\bold C)$.
\par

Denote by $B^4$ the open unit ball
$\{ (a,b) \in \bold C^2 \hskip.1cm \vert \hskip.1cm
\vert a \vert^2 + \vert b \vert^2 < 1 \}$,
by $\overline{B^4} = B^4 \sqcup \bold S^3$ its closure,
and by $\gamma$ the restriction to  $\overline{\bold B}^4$
of the inverse of the homeomorphism $g'$ defined above.
We modify $\gamma$ to obtain a mapping
$$
f \, : \, \overline{\bold B}^4 \hskip.2cm \longrightarrow \hskip.2cm
\Cal C^{\text{covol} \ge 1}(\bold C) \cup \Cal C_{\bold R}(\bold C)
$$
defined as follows:
\roster
\item"---"
$f(0,0) = \{0\}$,
\item"---"
if $(a,b) \in \bold S^3 \smallsetminus T$, then
$f(a,b) = \frac{1}{\sqrt{ \text{covol}(\gamma(a,b))}} \gamma(a,b)$,
\item"---"
if $(a,b) \in B^4 \smallsetminus \{(0,0)\}$
is in the $\bold R^*_+$--orbit of a point 
$(a_1,b_1) \in \left( \bold S^3 \smallsetminus T \right)$, then 
$f(a,b) = \frac{1}{\sqrt{h(a,b)}} \gamma(a_1,b_1)$
for an appropriate factor $h(a,b) \in ]0, \text{covol}(\gamma(a_1,b_1)) [$,
\item"---"
if $(a,b) \in T$ and $\gamma(a,b) = \bold Z \omega$,
then $f(a,b) = \bold R \omega$,
\item"---"
if $(a,b) \in B^4 \smallsetminus \{(0,0)\}$
is in  the $\bold R^*_+$--orbit of a point 
$(a_1,b_1) \in  T$, then
\newline
$f(a,b) = \frac{1}{\sqrt{h(a,b)}} \gamma(a_1,b_1)$
for an appropriate factor $h(a,b) \in ]0, \infty [$.
\endroster
For a precise definition of the function $h$,
see \cite{HuPo--79} and \cite{BrHK}.
It can be shown that $f$ is a homeomorphism.
\par

Moreover, $f$ can be extended to a mapping
$\bold C^2 \cup \{\infty\} \longrightarrow \Cal C(\bold C)$
by defining for $(a,b)$ outside $\bold B^4$
$$
f(a,b) \, = \, \big( f(\sigma(a,b) \big)^{\sharp}
$$
where $\sigma$ denotes the inversion 
$(a,b) \longmapsto \frac{(a,b)}{\vert a \vert^2 + \vert b \vert^2}$
of $\bold C^2$ fixing $\bold S^3$
and where the  {\it dual}
of a closed subgroup $C$ of $\bold C$
is defined by
$$
C^{\sharp} \, = \, \{ z \in \bold C \hskip.2cm \vert \hskip.2cm
\operatorname{Im}(\overline z \hskip.1cm c) \in \bold Z
\hskip.2cm \text{for all} \hskip.2cm c \in C \}    .
$$
(Observe that $C^{\sharp} = C$ if and only if $C$ 
is either a unimodular lattice 
or a subgroup isomorphic to $\bold R$.)
The one--point compactification $\bold C^2 \cup \{\infty\}$
of $\bold C^2$ can be identified to the $4$--sphere $\bold S^4$,
and we have: 

\proclaim{4.~Theorem \cite{HuPo--79}}
The mapping $f : \bold S^4 \longrightarrow \Cal C (\bold C)$
is a homeomorphism.
\endproclaim

\demo{Comments on Theorem~4}
By the homeomorphisms $f$, the equator $\bold S^3$ of $\bold S^4$
corresponds to the union 
of the set $\Cal L^{\text{unimod}}(\bold C)$ of unimodular lattices
and the set $\Cal C_{\bold R}(\bold C)$ of subgroups isomorphic to $\bold R$,
the latter corresponding to a trefoil knot $T \subset \bold S^3$.
\par
\Refs\nofrills{}
This cannot be seen using only the homeomorphism $g'$ 
of Subsection IV.1. Indeed, 
whenever a point $(a,b) \in \bold C^2 \smallsetminus \Sigma$
converges towards a point $(a_{\lim},b_{\lim}) \in \Sigma$,
the  corresponding  $(g')^{-1}(a,b)$ has a volume
which tends to infinite.
When  $(a,b)$ is rescaled in such a way that
the corresponding lattice is unimodular,
then $(a,b)$ escapes any compact subset of $\bold C^2$,
the minimal norm of the corresponding lattice tends to $0$,
and the lattice  itself 
tends inside $\Cal C (\bold C)$
to a subgroup isomorphic to $\bold R$.
\endRefs
\par
If we view $\bold S^4$ as the suspension of its equator $\bold S^3$, 
the two--dimensional sphere $\bold S^2$ of Subsection IV.1, which corresponds
to the complement of $\Cal L (\bold C)$ in $\Cal C (\bold C)$,
corresponds to the suspension of $T$.
A pole  has a typical neighbourhood inside  $\bold S^2$
which is a cone over the knot $T$
(other points of $\bold S^2$ have typical neighbourhoods in $\bold S^2$
which are cones over unknotted closed curves inside $\bold S^3$);
in particular, the embedding of this $\bold S^2$ in the total space $\bold S^4$
is not locally flat at the South and North Poles.
\par
We know that the space $\Cal L^{\text{unimod}}(\bold C)$
is not simply connected; indeed, it is an aspherical space with fundamental group
the inverse image of $SL_2(\bold Z)$ 
in the universal covering group of $SL_2(\bold R)$. 
An oriented closed curve $\ell$ in $\Cal L^{\text{unimod}}(\bold C)$
can be viewed as inside $\bold S^3$ and disjoint from the trefoil knot $T$
(oriented in some way),
so that there is an associated 
{\it linking number} $\operatorname{link}(\ell,T)$.
An interesting particular case is that of a periodic orbit of the
{\it geodesic flow} on $SL_2(\bold R) / SL_2(\bold Z)$,
viewed as the unit tangent bundle of the {\it modular surface}
$\bold D^2 / PSL_2(\bold Z)$,
as explained in \cite{Ghys--07}.
\enddemo

\bigskip
\head{\bf
V.~On the Heisenberg group
}\endhead

The Heisenberg group is a $3$--dimensional nilpotent Lie group,
which is connected and simply connected. 
We use the model
\footnote{
An element $(x+iy,t)$ corresponds to a matrix
$ \left( \matrix 1 & x & t + \frac{1}{2}xy \\ 0 & 1 & y \\ 0 & 0 & 1 \endmatrix \right)$
in the matrix picture
$\left( \matrix 1 & \bold R & \bold R \\ 0 & 1 & \bold R \\ 0 & 0 & 1 \endmatrix \right)$
of the Heisenberg group.
}
$$
H \, = \, \bold C \times \bold R
%\hskip.5cm \text{with product defined by} \hskip.5cm
\hskip.5cm \text{with} \hskip.5cm
(z_1,t_1)(z_2,t_2) \, = \, 
\Big( z_1+z_2 \hskip.1cm , \hskip.1cm
 t_1+t_2 + \frac{1}{2} \operatorname{Im}(z_1\overline{z_2}) \Big) ,
$$
we denote by $p : H \longrightarrow \bold C$
the projection $(z,t) \longmapsto z$,
and we identify $\bold R$ to $\{0\} \times \bold R \subset H$.
The following is a collection of easily verified properties.

\proclaim{5.~Proposition}
(i) The subset $\bold R$ of $H$ is both the commutator subgroup and the centre
$Z(H)$ of the Heisenberg group; moreover, commutators are given by
$$
[(z_1,t_1),(z_2,t_2)] \, = \,  (0,\operatorname{Im}(\overline{z_1}z_2)) .
$$
\par
(ii) Maximal abelian subgroups of $H$ are of the form
$\bold R z_0 \times R$, with $z_0 \in \bold C^*$.
\par
(iii) For any $C \in \Cal C (H)$ with $C \cap Z(H) \ne \emptyset$,
the projection $p(C)$ is a closed subgroup of~$\bold C$
(this applies in particular to non--abelian closed subgroups of $H$).
\par
(iv) For any non--abelian $C \in \Cal C (H)$,
either $Z(H) \subset C$ or $C \in \Cal L (H)$.
\par
(v) For any lattice $\Lambda$ in $H$,
the commutator subgroup $[\Lambda,\Lambda]$ 
is a subgroup of finite index in the centre $Z(\Lambda)$.
\par
(vi) The automorhism group of $H$ is a semi--direct product
$(H/Z(H)) \rtimes GL_2(\bold R)$, 
with $H/Z(H)$ the group of inner automorphisms 
and with $GL_2(\bold R)$ acting on $H$ by
$$
\left( \left( \matrix a & b \\ c & d \endmatrix \right)
\hskip.1cm , \hskip.1cm (x+iy,t) \right)
\hskip.3cm \longmapsto \hskip.3cm
\Big( (ax+by) + i(cx+dy) \hskip.1cm , \hskip.1cm (ad-bc)^2 t \Big) .
$$
\endproclaim
\bigskip

For any positive integer $n$, we denote by $\Cal L_n(H)$ the
subspace of lattices $\Lambda \in \Cal L (H)$ with 
$[\Lambda,\Lambda]$ of index $n$ in $Z(\Lambda)$.
We denote by $\Cal L_{\infty}(H)$ the subspace of $\Cal C (H)$
of subgroups of the form $p^{-1}(L)$, with $L \in \Cal L (\bold C)$.

\proclaim{6.~Examples} For any $n \ge 1$, the subgroup $\Lambda_n$ of $H$
generated by
$(1,0)$, $(i,0)$, and $(0,1/n)$ is in $\Cal L_n (H)$. Moreover,
$$
\Lambda_n \, = \, 
\bold Z [i] \times \frac{1}{n}  \bold Z 
$$
if $n$ is even and
$$
\Lambda_n \, = \, 
\left\{ \left(x+iy , \frac{t}{2n} \right) \in \bold Z [i] \times \frac{1}{2n} \bold Z
\hskip.2cm \Big\vert \hskip.2cm
xy \equiv t \pmod{2} \right\}
$$
if $n$ is odd. In the matrix picture for $H$, we have
$\Lambda_n = \left( 
\matrix 1 & \bold Z & \frac{1}{n}\bold Z \\ 0 & 1 & \bold Z \\ 0 & 0 & 1 \endmatrix
\right)$
for all $n$.
\par
The subgroup of $H$ generated by $(1,0)$, $(i,0)$, and $Z(H)$
is in $\Cal L_{\infty}(H)$.
\endproclaim

\bigskip
\head{\bf
VI.~Three subspaces of $\Cal C (H)$
}\endhead

The space $\Cal C (H)$ consists of several parts which,
as it is the case for $\Cal C (\bold C)$, are relatively easy to describe,
and their gluing is harder. From Proposition~5, 
the space $\Cal C (H)$ consists of three parts:
\roster
\item"$\to$"
the space $\Cal A (H)$ of closed abelian subgroups,
\item"$\to$"
the space $\cc (H)$ of closed subgroups of $H$ containing $Z(H)$, 
\item"$\to$" 
the space of lattices $\Cal L (H) = \bigsqcup_{n=1}^{\infty} \Cal L_n(H)$.
\endroster
Moreover,  $\Cal L (H)$ is disjoint from
$\Cal A (H) \cup \cc (H)$.

\medskip
\subhead
VI.1.~The space $\Cal A (H)$ of abelian subgroups
\endsubhead

From Proposition 5.ii, the space $\Cal C_{\bold R^2} (H)$
of maximal abelian subgroups of $H$
is homeomorphic to
the projective line $\Cal C_{\bold R} (\bold C) \approx \bold P^1$ 
of lines in the space $H/Z(H) \approx \bold R^2$.
For a closed abelian subgroup $A$ of $H$, there are two cases to distinguish:
\roster
\item"(i)" if $A \nsubseteq Z(H)$, then $A$ is contained 
in an unique maximal abelian subgroup  of $H$,
\item"(ii)" if $A \subseteq Z(H)$, then $A$ is contained 
in all maximal abelian subgroups of $H$.
\endroster
If one could identify coherently each maximal abelian subgroup of $H$ 
to a copy of $\bold R^2$, the space $\Cal A (H)$ would simply be the quotient of 
$\Cal C_{\bold R^2} (H) \times \Cal C (\bold C) 
\approx \bold S^1 \times \bold S^4$
by the relation identifying each of the circle
$(A_{\max},C)$ to a point 
(for $A_{\max}$ in the circle $\Cal C_{\bold R^2} (H)$,
and for  $C \subset Z(H)$  fixed).
This identification is possible locally, but not globally;
indeed, even once $H$ is furnished with some Riemannian structure, 
each $A_{\max}$ can be identified to $\bold C$ in two ways,
corresponding to the two orientations of $A_{\max}$.
However, the following can be shown (Proposition 6.1.i of \cite{BrHK}):

\proclaim{7.~Proposition}
There exists in $\bold S^4$ a tame closed interval $I$,
corresponding in $\Cal C (\bold C)$ to the subspace
of closed subgroups contained in $\bold R$,
such that
$$
\Cal A (H) \, \approx \,
\Big(S^1 \times S^4\Big) \big/ 
\Big( (\varphi , x)\sim (\varphi' , x),  
\hskip.2cm \text{for all $\varphi,\varphi' \in \bold S^1$ and $x \in I$} \Big).
$$
In particular, $\Cal A (H)$ is a space of dimension $5$.
\endproclaim

\medskip
\subhead
VI.2.~The space $\cc (H)$ of groups containing the centre $Z(H)$
\endsubhead

As a straightforward consequence of the results of Section~IV, we have:

\proclaim{8.~Proposition}
The space $\cc (H)$ is homeomorphic to $\Cal C (\bold C)$,
namely to a $4$--dimensional sphere. 
Moreover, the intersection $\Cal A (H) \cap \cc (H)$ is embedded in $\cc (H)$
as a closed $2$--disc. 
\endproclaim

More precisely, this intersection $\Cal A (H) \cap \cc (H)$ 
is the lower hemisphere of the $2$--sphere
described in the comments which follow Theorem~4.
Observe that the complement of this $2$--sphere in $\cc (H)$
is precisely the space $\Cal L_{\infty}(H)$ of subgroups of the form
$p^{-1}(L)$ for $L \in \Cal L (\bold C)$, as defined
just before Example~6.

\medskip
\subhead
VI.3.~The space $\Cal L (H)$ of lattices in $H$
\endsubhead

Consider an integer $n \ge 1$, a lattice $\Lambda_0 \in \Cal L_n (H)$,
and denote by $L_0 = p(\Lambda_0) \in \Cal L (\bold C)$ its projection.
Choose a positively oriented basis $(z_0,z'_0)$ of $L_0$
and two points $(z_0,t_0) \in \Lambda_0 \cap p^{-1}(z_0)$, 
$(z'_0,t'_0) \in \Lambda_0 \cap p^{-1}(z'_0)$.
The infinite cyclic subgroup $[\Lambda_0,\Lambda_0]$
is generated by the commutator
$$
[(z_0,t_0),(z'_0,t'_0)] \, = \,  (0,\operatorname{Im}(\overline{z_0}z'_0))
= (0,\operatorname{vol}(\bold C / L_0)) ,
$$
and $\Lambda_0$ itself is generated by $(z_0,t_0)$, $(z'_0,t'_0)$,
and $(0,\frac{1}{n}\operatorname{vol}(\bold C / L_0))$.
\par

For any nearby lattice $\Lambda \in \Cal L_n (H)$,
there exist a unique positively oriented basis $(z,z')$ 
of $L  \Doteq p(\Lambda)$ near $(z_0,z'_0)$
and numbers $t,t'$ near $t_0,t'_0$ such that
$(z,t)$, $(z',t')$,
and $(0,\frac{1}{n}\operatorname{vol}(\bold C / L))$
generate~$\Lambda$.
Moreover, any lattice in $\Cal L_n (H)$ with projection $L$
is generated by elements of the form 
$(z,t)$, $(z',t')$,
and $(0,\frac{1}{n}\operatorname{vol}(\bold C / L))$,
for some $t,t' \in [0,\frac{1}{n}\operatorname{vol}(\bold C / L)[$.
\par

It is easy to check that the subgroup $GL_2(\bold R)$
of $\operatorname{Aut}(H)$, see Proposition 5.vi,
operates on $\Cal L_n (H)$ in such a way that,
for any orbit $\Cal O$, we have
$$
\{ L \in \Cal L (\bold C) \hskip.1cm \vert \hskip.1cm
L = p(\Lambda) \hskip.2cm \text{for some} \hskip.2cm \Lambda \in \Cal O \}
\, = \, \Cal L (\bold C) .
$$
Moreover,  the group of inner automorphisms
operates transitively on the set of possible choices for $t,t'$.
We have essentially proved the first two claims of the following proposition;
for details and for the last claim, see \cite{BrHK, Section~7}.

\proclaim{9.~Proposition} For each $n \ge 1$, the space $\Cal L_n (H)$ is both
\par
(i) a torus bundle with base space $\Cal L (\bold C)$,
\par
(ii) a homogeneous space of the $6$--dimensional Lie group
$\operatorname{Aut}(H)$ by a discrete subgroup.
% \par
For $n,n' \ge 1$, the spaces $\Cal L_n (H), \Cal L_{n'} (H)$ 
are homeomorphic to each other; 
moreover, the torus bundles
$\Cal L_n (H) \longrightarrow \Cal L (\bold C)$ and
$\Cal L_{n'} (H) \longrightarrow \Cal L (\bold C)$ 
are isomorphic.
\endproclaim

\medskip
\subhead
VI.4.~Summing up for $\Cal A (H)$, $\cc (H)$, and $\Cal L (H)$
\endsubhead

Let $\left( \varphi_s \right)_{s > 0}$ be the one--parameter
subgroup of $\operatorname{Aut}(H)$ defined by
$\varphi_s (z,t) = (sz,s^2t)$. 
We have $\lim_{s \to \infty} \varphi_s (D) = \{e\}$
for any discrete subgroup $D$ of $H$ and 
% \footnote{
% The existence of such discrete subgroups
% converging to the whole group
% has been used to study some investigations
% of harmonic analysis \cite{Delm--06}.
% }
$\lim_{s \to 0} \varphi_s (\Lambda) = H$
for any $\Lambda \in \Cal L (H)$.
Since $\Cal A (H)$ and $\cc (H)$ are clearly arc--connected
(Propositions 7 and 8), and contain respectively $\{e\}$ and $H$,
it follows that $\Cal C (H)$ is arc--connected.
This is the very first part of the following theorem,
one of the two theorems which summarise the results in \cite{BrHK};
the space $\Cal C (H)$ is not locally connected
because any neighbourhood of a point in $\Cal L_{\infty} (H)$
is disconnected, containing points from $\Cal L_n (H)$
for $n$ large enough.

\proclaim{10.~Theorem (Theorem 1.3 in \cite{BrHK})}  
The compact space $\Cal C (H)$ 
is arc--connected but not locally connected.
It can be expressed as the union of 
the following three subspaces.
\roster
\item"(i)"
$\Cal L(H)$, which is open and dense in $\Cal C (H)$;
this has countably many connected components $\Cal L_n(H)$, 
each of which is homeomorphic to a fixed aspherical $6$--manifold that is a 
$2$--torus bundle over $\Cal L (\bold C) \approx GL_2(\bold R) / GL_2(\bold Z)$.
\item"(ii)"
$\Cal A(H)$, which is homeomorphic to  the space obtained from
$\bold S^4\times\bold P^1$ by fixing a tame arc $I\subset\bold S^4$
and collapsing each of the circles 
$\{\{i\} \times \bold P^1 : i\in I\}$
to a point.  
\item"(iii)"
$\cc(H)$, from which there
 is a natural homeomorphism to  $\bold S^4$;
the complement of $\Cal L_{\infty}(H)$ in $\cc (H)$ is a 
$2$--sphere $\Sigma^2\subset \bold S^4$
(which fails to be locally flat at  two points).
\endroster
The union $\Cal A(H) \cup \cc(H)$ is the complement of $\Cal L (H)$ in $\Cal C (H)$.
The intersection
$$
\Cal A(H) \cap \cc(H) =
\{ C \in \cc (H) \mid p(C) \subset
\Cal C_{\{0\},\bold Z,\bold R}(\bold C) \}
$$ 
is  a closed $2$--disc in $\Sigma^2$. 
The space $\Cal L (H) \cup \Cal L_{\infty}(H)$
is precisely $\{ C \in \Cal C(H) \mid p(C) \in \Cal L (\bold C) \}$.

\iffalse
The frontier of $\Cal L_n(H)$, which is independent of $n$, 
is the union of $\Cal A(H)$ and of the complement $\Sigma^2$
of $\Cal L_{\infty}(H)$ in $\cc (H)$.
\fi

\smallskip

$\cc (H)$ is a weak retract of $\Cal C (H)$:
there exists a continuous map 
$f : \Cal C (H) \longrightarrow \bold S^4$, 
constant on $\Cal A(H)$, such that 
\footnote{
Here, $\simeq$ denotes  homotopy equivalence.
}
$f \circ j \simeq {\operatorname{id}}_{\bold S^4}$, where
$j : \bold S^4 \longrightarrow \cc(H)$ is the homeomorphism of (iii).
In particular, $\pi_4(\Cal C (H))$ surjects onto $\bold Z$.

\smallskip

The subspace $\Cal N (H)$ of normal closed subgroups of $H$
is the union of $\cc (H)$ (which 
is homeomorphic to $\Cal C (\bold C)
\approx \bold S^4$)
and  the closed interval $\{ C \in \Cal C (H) \mid C \subset Z(H) \}$,
attached to the sphere $\cc (H)$ by one of its endpoints.
\endproclaim 

We would like to know more generally when
$\Cal L (G)$ is dense in $\Cal C (G)$, 
say for a unimodular connected Lie group $G$.

% remarque de Derighetti, $\Cal L$ jamais dense dans $\Cal C$
% pour $G$ de Lie semi--simple ????

% ???? Kazhdan--Margulis ????

\bigskip
\head{\bf
VII.~The space  $\Cal C (H)$
}\endhead

As an illustration of the density of $\Cal L (H)$ in $\Cal C (H)$,
let us describe the following simple

\proclaim{11.~Example} Let $n \ge 1$ be fixed.
For any integer $k \ge 1$, let $A_k$ 
denote the subgroup of $\bold R \times \bold R$
generated by $(1,0)$ and $(-\frac{1}{k},1)$;
let $A$ denote the subgroup $\bold Z^2$ of $\bold R \times \bold R$;
and let $\Lambda_k$ denote the subgroup of $H$
generated by $A_k$ and $(-ik^2 n,0)$.
\par
Then $\Lambda_k \in \Cal L_n (H)$ for all $k \ge 1$ and 
$\lim_{k \to \infty} \Lambda_k = \lim_{k \to \infty} A_k =  A$ 
in $\Cal C (H)$.
\endproclaim

\demo{Remark, and consequence of this example} 
Here, $\bold R \times \bold R$ is viewed as a subgroup of
$\bold C \times \bold R = H$;
in particular, it should not be confused with $\bold C$.
(Inside $H = \bold C \times \bold R$, note that $\bold C$
is not a subgroup.)
Observe that $p(A_k)$ is a lattice in $\bold C$ for each $k \ge 1$,
but the projection of the limit, $p(A)$, is a cyclic subgroup of $\bold C$.
\par
It follows from Example~11 that the frontier of $\Cal L_n (H)$ in $\Cal C (H)$
contains the closure of the $\operatorname{Aut}(H)$--orbit of $A$,
and it is a fact
that this closure coincides with $\Cal A (H)$
(Proposition 6.1.ii in \cite{BrHK}).
Moreover, it can be seen
that the frontier of $\Cal L_n (H)$ contains $\cc (H)$.
\enddemo

\medskip

The previous argument shows part of the following theorem,
again copied from \cite{BrHK}, by which we will
end this account. Recall that
we denote here by $\Sigma^2$ the topologically embedded
sphere in $\cc (H)$ which corresponds to the
$\bold S^2 \subset \bold S^4$ in the comments on Theorem~4.
The symbols $\bold P^2$ and $\bold K$ stand respectively
for a real projective plane and a Klein bottle.

\proclaim{12.~Theorem (Theorem 1.4 in \cite{BrHK})}  
The spaces $\Cal L_n(H)$ are homeomorphic  
to a common aspherical homogeneous space, namely the quotient
of the  $6$--dimensional automorphism group
$\operatorname{Aut}(H)  \cong
\bold R^2  \rtimes  GL_2(\bold R)$
by the discrete subgroup $\bold Z^2\rtimes GL_2(\bold Z)$.
\par

The frontier of $\Cal L_n(H)$, which is independent of $n$,
consists of the following subspaces:
\roster
\item"(i)"
the trivial group $\{e\}$;
\item"(ii)"
$\Cal C_{\bold R}(H) \approx \bold P^2$;  
\item"(iii)"
$\Cal C_{\bold Z}(H) \approx \bold P^2 \times ]0,\infty[$;
\item"(iv)"
$\Cal C_{\bold R^2}(H) \approx \bold P^1$;
\item"(v)"
$\Cal C_{\bold R \oplus \bold Z}(H) \approx \bold K \times ]0,\infty[$,
which is a $(\bold P^1 \times ]0,\infty[)$--bundle over $\bold P^1$;
\item"(vi)"
$\Cal C_{\bold Z^2}(H)$,
which is a $(\bold S^4 \smallsetminus \Sigma^2)$--bundle over $\bold P^1$;
\item"(vii)"
$p_*^{-1}\left( \Cal C_{\bold R \oplus \bold Z}(\bold C) \right)$;
\item"(viii)"
the full group $H$.
\endroster
In particular, the frontier of $\Cal L_n(H)$ is the union of  $\Cal A(H)$
and the complement $\Sigma^2$ of $\Cal L_{\infty}(H)$ in $\cc (H)$;
the part $\Cal A (H)$ is itself  the union of the subspaces {\rm (i)} to {\rm (vi)},
and $\Sigma^2 \smallsetminus (\Sigma^2 \cap \Cal A (H))$ 
is itself the union of the subspaces {\rm (vii)} and {\rm (viii)}.
The frontier of $\bigcup_{n=1}^{\infty} \Cal L_n(H)$ further contains
\roster
\item"(ix)"  
$\Cal L_{\infty}(H)$.
\endroster
Each of these spaces, except~{\rm (vi)}, consists of finitely many
$\operatorname{Aut}(H)$--orbits.
\endproclaim

Observe that, as $\Cal L (H)$ is open dense, 
the spaces (i) to (ix) of Theorem~1.4
together with the spaces $\Cal L_n (H)$ for $n \ge 1$
constitute a partition of $\Cal C (H)$.

% \vskip1cm \centerline{??? Picture ???}

\enddocument